\documentclass[10pt,one side,emlines]{amsart}
\usepackage{amssymb,latexsym,xy,eucal,mathrsfs}
\usepackage[utf8]{inputenc}
\usepackage{tikz}
\usetikzlibrary{arrows.meta}
\usepackage{graphicx, graphics}
\usepackage{caption}
\textwidth=14cm \textheight=22cm \theoremstyle{plain}
\newtheorem{lem}{Lemma}[section]

\newtheorem{cnj}[lem]{Conjecture}
\theoremstyle{definition}

\begin{document}
	\baselineskip 15truept
	
	%\thanks{ This research is supported by Board of College and University Development, University of Pune, via the project SC-66}
	
	\title{Counter-example to Conjectures on Complemented Zero-Divisor Graphs of Semigroups}
	\maketitle
	\markboth{Anagha Khiste, Ganesh Tarte and Vinayak Joshi}{Counter-example to Conjectures}\begin{center}\begin{large} $^\text{a}$Anagha Khiste, $^\text{b}$Ganesh Tarte and $^\text{c}$Vinayak Joshi \end{large}\begin{small}\vskip.1in$^\text{a}$\emph{Department of Applied Science and Humanities, Indian Institute of Information Technology, Pune - 410507 }\\$^\text{b}$\emph{Department of Applied Sciences and Humanities, Pimpri Chinchwad College of Engineering, Pune - 411044}\\$^\text{c}$\emph{Department of Mathematics, Savitribai Phule Pune University, Pune - 411007, Maharashtra, India\\E-mail: avanikhiste@gmail.com, ganesh.tarte@pccoepune.org, vvjoshi@unipune.ac.in, vinayakjoshi111@yahoo.com }\end{small}\end{center}\vskip.2in

	\begin{abstract}
		In this paper, we are motivated by two conjectures proposed by C. Bender et al.\ in 2024, which have remained open questions. The first conjecture states that if the complemented zero-divisor graph \( G(S) \) of a commutative semigroup \( S \) with a zero element has the clique number three or greater, then the reduced graph \( G_r(S) \) is isomorphic to the graph \( G(\mathcal{P}(n)) \). The second conjecture asserts that if \( G(S) \) is a complemented zero-divisor graph with the clique number three or greater, then \( G(S) \) is uniquely complemented.
		
		In this work, we construct a commutative semigroup \( S \) with a zero element that serves as a counter-example to both conjectures. 
	\end{abstract}	

	\maketitle \noindent{ \small \textbf{Keywords}: Clique number, commutative semi-group,  complemented graph, reduced graph, uniquely complemented graph, zero-divisor graph.} 
    \maketitle
	
	\section{Introduction}
	The zero-divisor graph associated with a commutative ring was first introduced by I.~Beck~\cite{B}, with the goal of studying graph coloring. Later, Anderson and Livingston~\cite{AL} modified the definition so that the vertex set consists only of the nonzero zero-divisors of the ring. This construction was further extended to the context of semigroups and posets in~\cite{DD,DMS, JA 2}. 
	
In a graph $G$, two distinct vertices $a$ and $b$ are called \textit{orthogonal}, denoted $a\perp b$, if $a$ and $b$ are adjacent and there does not exist a vertex $c$ adjacent to both $a$ and $b$. Equivalently, adjacent vertices $a$ and $b$ are \textit{orthogonal} if the edge $a - b$ is not an edge of any triangle (3-cycle). A graph $G$ is called \textit{complemented} if for every vertex $a$, there exists a vertex $b$ such that $a\perp b$. A graph $G$ is called \textit{uniquely complemented} if $G$ is complemented and has the property that if $a, b, c$ are distinct vertices with $a\perp b$ and $a\perp c$, then $N(b) = N(c)$. This graph-theoretic property has been studied in the context of the zero-divisor graph associated with a commutative ring; \textit{see} \cite{L}. 
	
		Throughout the paper, let $S$ denote a finite commutative semigroup with $0$. The set of zero-divisors of $S$, denoted $Z(S)$, forms an ideal in $S$. The zero-divisor graph associated with $S$, denoted $G(S)$, is the graph whose vertices are given by the nonzero zero-divisors and  two distinct vertices $a$ and $b$ are	adjacent  when $ab = 0$.
	  A \textit{clique}	in a graph $G$ is a sub-graph of $G$ where every pair of distinct vertices in the sub-graph are adjacent. The	\textit{clique number} of a graph $G$ is the number of vertices in a maximal clique in $G$ and is denoted $\omega(G)$.
	
	Define an equivalence relation $\sim$ on the vertex set $V(G)$ of a graph $G$ by $a\sim b$ if	and only if $N(a) = N(b)$. Denote the equivalence class of $a\in V(G)$ by $[a]$. The \textit{reduced graph} of $G$, denoted $G_r$, is the graph whose vertex set is the set of equivalence classes under relation $\sim$, where two distinct vertices $[a]$, $[b]$ in the vertex set of $G_r$ are adjacent if and only if vertices $a$ and $b$ are	adjacent in the graph $G$; \textit{see} Bender et al.\cite{C}.

	Let $\mathcal{P}(n)$ denote the semigroup consisting of the power set of an $n$-element set, equipped with the operation of intersection.
	
	  In \cite{C}, Bender et al.  mentioned that	every known example of a zero-divisor graph $G(S)$, for finite $S$, which is complemented and has the clique number $n \geq 3$ has reduced graph isomorphic to $G(\mathcal{P}(n))$ and is also a uniquely complemented graph. The case where the order of $S$ is infinite is open. Hence they raised the following conjectures.

	\begin{cnj} If $G(S)$ is a complemented zero-divisor graph with the clique number $n \geq 3$, then $G(S)$ has the reduced graph $G_r$ which is isomorphic to  $G(\mathcal{P}(n))$.
		\end{cnj}
		\begin{cnj} If $G(S)$ is a complemented zero-divisor graph with the clique number $3$ or greater, then $G(S)$ is uniquely complemented.
			\end{cnj}
	To settle these conjectures, we explored the underlying structures of complemented zero-divisor graphs in commutative semigroups with zero.

	\section{Counter-example}
	Consider the set $S=\mathcal{P}(3)\times \mathbb{Z}_4=\{(X,\overline{a})~|~ X \in \mathcal{P}(3),~\overline{a} \in \mathbb{Z}_4 \}$. 
	Clearly, $S$ is a nonempty set, as $(\emptyset,\overline{0}) \in S$.  We define a binary operation, denoted by $``*"$ on the set $S$ such that $(X,\overline{a}) * (Y,\overline{b})=(X \cap Y,~~\overline{a}\times_4\overline{b})$, for $(X,\overline{a}), (Y,\overline{b}) \in S$. For simplicity, we denote the binary operation $\times_4$ on $\mathbb{Z}_4$ by $``."$. 

Since $(\mathcal{P}(3),\cap)$~ and ~ $ (\mathbb{Z}_4,.)$ are both commutative semigroups with zero elements, we can easily prove that $(S,*)$ is a commutative semigroup with zero element $(\emptyset, \overline{0})$.
	
Observe that,
$V(G(S))=\{(X,\overline{a}) \in S \setminus (\emptyset,\overline{0}) ~|~ (X,\overline{a})*(Y,\overline{b})=(\emptyset,\overline{0})$, \linebreak ~~\text{for some}~~ $(Y,\overline{b}) \in S \setminus\{(\emptyset,\overline{0})\}\}$. Clearly, $V(G(S)) \neq \emptyset$, since $(\emptyset, \overline{1})*(\{123\}, \overline{0})=(\emptyset, \overline{0})$. Further, we can easily prove that $V(G(S))=S\setminus \{(\emptyset, \overline{0}), (\{123\},\overline{1}), (\{123\},\overline{3})\}$. 

First, we claim that the zero-divisor graph $G(S)$ is complemented.	

Let $(X,\overline{a}) \in V(G(S))$ be any element. Thus $(X,\overline{a})\neq (\emptyset, \overline{0})$ and there exists $(Y,\overline{b}) \neq (\emptyset, \overline{0})$ such that $(X,\overline{a})*(Y,\overline{b})= (\emptyset, \overline{0}).$ Therefore, $X \cap Y=\emptyset$ and $\overline{a}.\overline{b}=\overline{0}$.
Since $(X,\overline{a}) \neq (\emptyset,\overline{0})$, we have $X \neq \emptyset ~~ \text{or} ~~ \overline{a} \neq \overline{0}$.
Thus, we have the following cases:

	\begin{enumerate} \item[Case (I) :] Let $X =\emptyset$ and $\overline{a} \in \{\overline{1},\overline{2},\overline{3}\}$. Then $(X,\overline{a})\in \{(\emptyset,\overline{1}), (\emptyset,\overline{2}), (\emptyset,\overline{3})\}$.\\ It is easy to observe that $(X,\overline{a})\perp (\{1,2,3\},\overline{0})$. 
		\item[Case (II) :] Let $X \neq \emptyset$ and $\overline{a}=\{\overline{0}\}$. 
	Clearly, one can see that $(X,\overline{0}) \perp (Y, \overline{b})$ where \linebreak $(Y, \overline{b}) \in \biggl\{ (\{1,2,3\}\setminus X,\; \overline{1}), (\{1,2,3\}\setminus X,\; \overline{2}), (\{1,2,3\}\setminus X, \;\overline{3})\biggr\}$.
		\item[Case (III) :] Let $X \neq \emptyset$ and $\overline{a} \in \{\overline{1},\overline{2},\overline{3}\}$.  Then we have $ (Y, \overline{b})\perp (\{1,2,3\}\setminus X,\overline{0})$, where $(Y, \overline{b}) \in\biggl\{(X,\overline{1}), (X,\overline{3})\biggr\}$ provided $\{1,2,3\}\setminus X\neq \emptyset$. 	Further, $(X,\overline{2}) \perp (Y, \overline{b})$ where $(Y, \overline{b}) \in\biggl\{ (\{1,2,3\}\setminus  X,\overline{0}), (\{1,2,3\}\setminus X,\overline{2})\biggr\}$.
	\end{enumerate} 
Thus in any case, for any element $(X,\overline{a})\in V(G(S))$ has a complement in $G(S)$.	Therefore, $G(S)$ is a complemented graph. Clearly, $G(S)$ contains a cycle of length $3$, namely, $(\{1\},\overline{0})-(\{2\},\overline{0})-(\{3\},\overline{0})-(\{1\},\overline{0})$. This implies that {\textbf{$\omega{(G(S))} \geq 3$}}. 
	
\vskip10pt	
\noindent\textbf{Proof of Conjecture 1.1:}

	Now, for the reduced graph $G_r(S)$, first we determine neighborhood $N((X,\overline{a}))$, for every vertex $(X,\overline{a})\in V(G(S))$ as follows; 
	
	\begin{enumerate}
		\item[A)] Let $X=\emptyset$ and $\overline{a}\in \{\overline{1},\overline{2}, \overline{3}\}$. \begin{enumerate}\item $N((\emptyset,\overline{1}))=\{(X,\overline{0})~|~X(\neq \emptyset) \in \mathcal{P}(3)\}$ \\ $~~~~=\{(\{1\},\overline{0}),(\{2\},\overline{0}),(\{3\},\overline{0}),(\{1,2\},\overline{0}),(\{1,3\},\overline{0})\},(\{2,3\},\overline{0}),(\{1,2,3\},\overline{0})$
		$=N((\emptyset,\overline{3}))$.
		
		\item $N((\emptyset,\overline{2}))=\{(X,\overline{0}),(X,\overline{2})~|~X(\neq \emptyset) \in \mathcal{P}(3)\}$ \\ $~~~~=\{(\{1\},\overline{0}),(\{2\},\overline{0}),(\{3\},\overline{0}),(\{1,2\},\overline{0}),(\{1,3\},\overline{0})\},(\{2,3\},\overline{0}),(\{1,2,3\},\overline{0}), \\ ~~~~~~~~~~~~~~ (\{1\},\overline{2}),(\{2\},\overline{2}),(\{3\},\overline{2}),(\{1,2\},\overline{2}),(\{1,3\},\overline{2})\},(\{2,3\},\overline{2}),(\{1,2,3\},\overline{2})$.
		\end{enumerate}
		\item[B)] Let $X\neq \emptyset$ and $\overline{a}=\overline{0}$. \\Thus, $N((X,\overline{0}))=\{(\emptyset,\overline{a}),(Y,\overline{b})~|~\overline{a} \in \{\overline{1},\overline{2},\overline{3}\}, \overline{b} \in \mathbb{Z}_4, Y(\neq \emptyset) \in \mathcal{P}(3)$, such that $X \cap Y=\emptyset \}$. 
		\begin{enumerate}
			\item \mbox{ $N((\{1\},\overline{0}))=\{(\emptyset,\overline{a}),(\{2\},\overline{b}),(\{3\},\overline{b}), (\{2,3\},\overline{b})~|~ \overline{a}\in \{\overline{1},\overline{2}, \overline{3}\} ~\text{and}~ \overline{b}\in \mathbb{Z}_4\}$.}
			\item \mbox{$N((\{2\},\overline{0}))=\{(\emptyset,\overline{a}),(\{1\},\overline{b}),(\{3\},\overline{b}), (\{1,3\},\overline{b})~|~ \overline{a}\in \{\overline{1},\overline{2}, \overline{3}\} ~\text{and}~ \overline{b}\in \mathbb{Z}_4\}$.}
			\item \mbox{$N((\{3\},\overline{0}))=\{(\emptyset,\overline{a}),(\{1\},\overline{b}),(\{2\},\overline{b}), (\{1,2\},\overline{b})~|~ \overline{a}\in \{\overline{1},\overline{2}, \overline{3}\} ~\text{and}~ \overline{b}\in \mathbb{Z}_4\}$.}
			\item \mbox{$N((\{1,2\},\overline{0}))=\{(\emptyset,\overline{a}),(\{3\},\overline{b})~|~ \overline{a}\in \{\overline{1},\overline{2}, \overline{3}\} ~\text{and}~ \overline{b}\in \mathbb{Z}_4\}$.}
			\item \mbox{$N((\{1,3\},\overline{0}))=\{(\emptyset,\overline{a}),(\{2\},\overline{b})~|~ \overline{a}\in \{\overline{1},\overline{2}, \overline{3}\} ~\text{and}~ \overline{b}\in \mathbb{Z}_4\}$.}
			\item \mbox{$N((\{2,3\},\overline{0}))=\{(\emptyset,\overline{a}),(\{1\},\overline{b})~|~ \overline{a}\in \{\overline{1},\overline{2}, \overline{3}\} ~\text{and}~ \overline{b}\in \mathbb{Z}_4\}$.}
			\item \mbox{$N((\{1,2,3\},\overline{0}))=\{(\emptyset,\overline{a})~|~ \overline{a}\in \{\overline{1},\overline{2}, \overline{3}\}\}$.}
		\end{enumerate}
		Clearly, all the neighborhoods $N((\{1\},\overline{0}))$, $ N((\{2\},\overline{0}))$, $ N((\{3\},\overline{0}))$, $ N((\{1,2\},\overline{0}))$, $ N((\{2,3\},\overline{0}))$, $ N((\{1,3\},\overline{0}))$ are distinct.
		 
		\item[C)] Let $X (\neq \emptyset)$, and $\overline{a}\in \{\overline{1}, \overline{3}\}$.\\ Thus, $N((X,\overline{a}))=\{(Y,\overline{0})~|~ Y (\neq \emptyset) \in \mathcal{P}(3) ~\text{such that}~ X \cap Y=\emptyset\}$. 
		\begin{enumerate}
			\item	$N((\{1\},\overline{a}))=\{(\{2\},\overline{0}),(\{3\},\overline{0}),(\{2,3\},\overline{0})\}$
			\item	$N((\{2\},\overline{a}))=\{(\{1\},\overline{0}),(\{3\},\overline{0}),(\{1,3\},\overline{0})\}$
			\item	$N((\{3\},\overline{a}))=\{(\{1\},\overline{0}),(\{2\},\overline{0}),(\{1,2\},\overline{0})\}$
			\item	$N((\{1,2\},\overline{a}))=\{(\{3\},\overline{0})\}$
			\item	$N((\{1,3\},\overline{a}))=\{(\{2\},\overline{0})\}$
			\item	$N((\{2,3\},\overline{a}))=\{(\{1\},\overline{0})\}$	
		\end{enumerate}
		
		Clearly, all the neighborhoods $N((\{1\},\overline{a})) $, $  N((\{2\},\overline{a}))$, $  N((\{3\},\overline{a}))$, $  N((\{1,2\},\overline{a})) $, $  N((\{1,3\},\overline{a})) $, $  N((\{2,3\},\overline{a}))$ are distinct. 
		
		\item[D)] Let $X (\neq \emptyset)$, $\overline{a}=\overline{2}$. \\Thus, $N((X,\overline{2}))=\{(\emptyset,\overline{2}),(Y,\overline{a})~|~ Y (\neq \emptyset) \in \mathcal{P}(3)~\text{ such that}~ X \cap Y=\emptyset~\text{and}~\overline{a}\in \{\overline{0}, \overline{2}\}\}$. 
		\begin{enumerate}
			\item	$N((\{1\},\overline{2}))=\{(\emptyset,\overline{2}),(\{2\},\overline{a}),(\{3\},\overline{a}),(\{2,3\},\overline{a})~|~\overline{a}\in \{\overline{0}, \overline{2}\}\}$
			\item $N((\{2\},\overline{2}))=\{(\emptyset,\overline{2}),(\{1\},\overline{a}),(\{3\},\overline{a}),(\{1,3\},\overline{a})~|~\overline{a}\in \{\overline{0}, \overline{2}\}\}$
			\item $N((\{3\},\overline{2}))=\{(\emptyset,\overline{2}),(\{1\},\overline{a}),(\{2\},\overline{a}),(\{1,2\},\overline{a})~|~\overline{a}\in \{\overline{0}, \overline{2}\}\}$
			\item	$N((\{1,2\},\overline{2}))=\{(\emptyset,\overline{2}),(\{3\},\overline{a})~|~\overline{a}\in \{\overline{0}, \overline{2}\}\}$
			\item	$N((\{1,3\},\overline{2}))=\{(\emptyset,\overline{2}),(\{2\},\overline{a})~|~\overline{a}\in \{\overline{0}, \overline{2}\}\}$
			\item	$N((\{2,3\},\overline{2}))=\{(\emptyset,\overline{2}),(\{1\},\overline{a})~|~\overline{a}\in \{\overline{0}, \overline{2}\}\}$
			\item	$N((\{1,2,3\},\overline{2}))=\{(\emptyset,\overline{2})\}$
			\end{enumerate}
	
	Clearly, all the neighborhoods $N((\{1\},\overline{2}))$, $ N((\{2\},\overline{2})) $, $ N((\{3\},\overline{2}))$, $ N((\{1,2\},\overline{2})) $, $ N((\{1,3\},\overline{2}))$, $ N((\{2,3\},\overline{2}))$, $ N((\{1,2,3\},\overline{2}))$ are distinct. 
	\end{enumerate}
	
Hence, $[(X, \overline{1})]=[(X, \overline{3})]$. 

\noindent\mbox{Thus $V(G_r)=\biggl\{[(\emptyset,\overline{2})], [(X,\overline{a})], [(Y,\overline{1})]~|~X(\neq \emptyset)\in \mathcal{P}(3), Y\in \mathcal{P}(3)\setminus \{1,2,3\} ~\text{and}~\overline{a}\in \{\overline{0},\overline{2}\}\biggr\}$.
}
 Therefore, $|V(G_r(S))|=22$.

	Now, if we assume that, the reduced graph $G_r(S)$ is isomorphic to the graph $G(\mathcal{P}(n))$, then we have $22=|V(G_r(S))|=|V(G(\mathcal{P}(n)))|=2^n-2$. This implies that $2^n=24$, which is contradiction to the fact that $n \in \mathbb{N}$. \\ Therefore, $G_r (S)$ is not isomorphic to the graph $G(\mathcal{P}(n))$, for any $n \in \mathbb{N}$. 
	
	Thus Conjuncture 1.1 fails.  
	
\vskip10pt
\noindent \textbf{Proof of Conjecture 1.2:}

	We now claim that $G(S)$ is not a uniquely complemented graph.
	
	Let $u = (\emptyset, \overline{2})$, $v = (\{1,2,3\}, \overline{0})$, and $w = (\{1,2,3\}, \overline{2})$ be elements of $V(G(S))$. Observe that  
$	u * v = (\emptyset, \overline{0}),$
	and hence $u$ is adjacent to $v$ in $G(S)$. Suppose there exists an element $t = (X, \overline{a}) \in V(G(S))$, with $t \neq u$ and $t \neq v$, such that $t$ is adjacent to both $u$ and $v$, i.e.,  
	$t * u = t * v = (\emptyset, \overline{0}).$
	This implies $X = \emptyset$ and $\overline{a} \in \{\overline{0}, \overline{2} \}$.
	
	If $\overline{a} = \overline{0}$, then $t = (\emptyset, \overline{0})$, which is not a vertex of $G(S)$, since this element is the zero element of the semigroup and thus excluded from the vertex set.
	
	If $\overline{a} = \overline{2}$, then $t = (\emptyset, \overline{2}) = u$, contradicting the assumption $t \neq u$.
	
	Therefore, no such $t$ exists, and the edge $u$--$v$ is not part of a 3-cycle. Thus, $u \perp v$ in $G(S)$. Similarly, it can be shown that $u \perp w$ in $G(S)$.
	
	Next, we examine the neighborhoods of $v$ and $w$. Clearly,
	$N(v) = \{ (\emptyset, \overline{1}), (\emptyset, \overline{2}), (\emptyset, \overline{3}) \},
$ \; $
	N(w) =  \{ (\emptyset, \overline{2}) \}.
$
	Since $N(v) \neq N(w)$, it follows that $G(S)$ is not a uniquely complemented graph.  Therefore, Conjecture 1.2 fails.

\end{document}